\documentclass[letterpaper, 10 pt, conference]{ieeeconf}  

\IEEEoverridecommandlockouts                              
 \usepackage{amsmath}
\usepackage{amssymb}
\usepackage{amsfonts}
\usepackage{stfloats}
\usepackage{bbm}
\usepackage{color}
\usepackage{enumerate}	
\usepackage{wrapfig}
\usepackage{graphicx,epsfig}
\usepackage{bm}
\usepackage[margin=1in,papersize={8.5in,11in}]{geometry}
\usepackage{times}
\usepackage{multirow}
\usepackage{color}
\usepackage{xcolor}
\usepackage{subcaption}

\title{\LARGE \bf
Modeling and Control of Large-Scale Adversarial Swarm Engagements
}

\author{Theodoros Tsatsanifos,
        Abram H. Clark,
        Claire Walton,
         Isaac Kaminer,
        and~Qi Gong%
\thanks{Theodoros Tsatsanifos, Isaac Kaminer, and Claire Walton, are with the Department of Mechanical and Aerospace Engineering,
		Naval Postgraduate School, Monterey, CA, 93943 USA (e-mail:~theodoros.tsatsanifos.gr@nps.edu;~kaminer@nps.edu;~clwalton1@nps.edu).}
\thanks{Abe Clark is with the Department of Physics, Naval Postgraduate School (e-mail:~abe.clark@nps.edu).}
\thanks{Qi Gong is with the Department of Applied Mathematics, University of California Santa Cruz, Santa Cruz, CA 95064 USA (e-mail:~qgong@ucsc.edu).}
\thanks{Manuscript received on June 15, 2020.}
}

\begin{document}

\maketitle
\thispagestyle{empty}
\pagestyle{empty}

\begin{abstract}
We theoretically and numerically study the problem of optimal control of large-scale autonomous systems under explicitly adversarial conditions, including probabilistic destruction of agents during the simulation. Large-scale autonomous systems often include an adversarial component, where different agents or groups of agents explicitly compete with one another. An important component of these systems that is not included in current theory or modeling frameworks is random destruction of agents in time. In this case, the modeling and optimal control framework should consider the attrition of agents as well as their position. We propose and test three numerical modeling schemes, where survival probabilities of all agents are smoothly and continuously decreased in time, based on the relative positions of all agents during the simulation. In particular, we apply these schemes to the case of agents defending a high-value unit from an attacking swarm. We show that these models can be successfully used to model this situation, provided that attrition and spatial dynamics are coupled. Our results have relevance to an entire class of adversarial autonomy situations, where the positions of agents and their survival probabilities are both important.

\end{abstract}
\section{Introduction}  

Rapid technological advances have made large-scale networked swarms of autonomous agents a reality. Questions about robustness and resilience are relevant in the context of designing algorithms to control individual agents~\cite{evers2014robust} or collective swarms~\cite{reynolds1987flocks,vicsek,leonard2001virtual,brambilla2013swarm}. In many situations involving multiple distinct autonomous systems, there is a natural adversarial component to their operational environment: agents are (1) trying to accomplish a task while (2) minimizing their risk of crashing or otherwise being neutralized. For example, Unmanned Aerial Vehicles (UAVs) delivering packages may have to operate in a crowded airspace and, therefore, the onboard control algorithms should guarantee efficiency of the delivery as well as for mitigate for the risk of crashing (into another UAV or some other object). These interactions may be indirectly adversarial (optimal performance for one swarm might necessitate sub-optimal performance for another) or directly adversarial (competing goals). A recent review~\cite{TR2018} discusses adversarial control and gives bio-inspired examples such as birds of prey herding bird flocks \cite{herdingbirds}, dolphins hunting \cite{Egerstedt}, and sheep-dogs \cite{sheepdogs1, sheepdogs2}. Herding strategies such as these use the swarm's own response strategy against it, to generate goal outcomes such as containment. 

However, these studies do not include any agent attrition, which would be relevant for collisions or if the agents have means of neutralizing other agents \cite{walton2018optimal}.  An adversarial situation in which swarm members can be removed creates additional problems with changing network topology, swarm size, and intra-swarm dynamics. The attrition of agents becomes coupled with spatial motion of the agents, as the removal of a swarm member changes both intra-swarm dynamics as well as potentially interactions between agents in different swarms. UAVs operating in adversarial environments must employ algorithms that maximize some measure of success while accounting for their survival probabilities and the changing intra-swarm dynamics caused by agent removal. This latter piece—the coupling of survival probabilities with swarm dynamics—is not treated by any current theoretical framework.

In this paper, we develop a new theoretical framework for modeling and control of large-scale adversarial autonomy involving both agent dynamics and agent attrition. We consider a situation involving a set of N agents trying to accomplish some task. Depending on their position, these agents have some rate of attrition, meaning that the index set of active agents at each instant of time is reduced during the task execution in a way that is inherently stochastic. We then pose a novel optimal control problem, where we seek to maximize the probability of completing a task under these conditions.  We test this framework using agent-based numerical simulations in a swarm-versus-swarm adversarial engagement. 
Our results demonstrate a way forward for using direct methods of optimal control to solve a class of adversarial autonomy problems where agent attrition is of critical importance.


\section{Modeling and Optimization Framework for Adversarial Swarming} \label{sec:sec_3.1}
\label{sec:modeling} 

For true autonomous systems (i.e., no human operator), each agent typically operates in a deterministic fashion. However, attrition is much more likely to be random or probabilistic. For swarms or other large-scale networked autonomous systems, the agent dynamics may depend on the behavior of the other agents, especially those nearby. As time advances, agents may be killed randomly, which will affect the dynamics of other agents, causing a ripple effect on the swarm. Thus random attrition introduces some inherently stochastic features to the global behavior.

To capture the varying numbers of agents, we consider an index set $I(t_k)$ that labels agents survived at time instance $t_k$. We will study the specific example of $M$ defending agents protecting an asset against an attacking swarm of $N$ agents. At the initial time $t_0$, $I(t_0)=\{1,2,\ldots, N,N+1,\ldots, N+M\}$ to include all attackers and defenders at the start of the engagement. 
Clearly, for all $t_k>t_0$, $I(t_k)\subseteq I(t_0)$ due to the attrition. The discrete dynamics of all the agents in the swarm on swarm scenario can be summarized by
\begin{eqnarray}
 z(t_{k+1}) & = & \phi^k( z(t_k), u(t_k), I(t_k)), 
\label{eq:agent_dyn}
\end{eqnarray}
where $z(t_{k+1})$ and the corresponding controls are aggregated states/controls for agents in the index set at time $t_k$:
\[
z(t_{k+1}) = \bigcup_{i\in I(t_{k})} z_i
\]
and $\phi^k$ are the collective dynamics only for agents in the indicator set $ I(t_k)$ .

The superscript $k$ in $\phi^k$ emphasizes the time dependence of the dynamics function, which changes with each change in $I(t_k)$. More than just time dependency, it emphasizes the changing {\em dimension} of the state vector and its dynamics function as the index set changes.

Swarm dynamics (\ref{eq:agent_dyn}) appear to be deterministic. However, over the entire swarm engagement, the dynamics $\phi^k$ can change in a random fashion depending, for example, on the probability of survival of each agent. To capture such stochastic behavior, we introduce the following dynamical update of the index set.
\begin{eqnarray}
I(t_{k+1}) & = & \psi (I(t_k), z(t_k),\omega(t_k)), \label{eq:I_dyn}
\end{eqnarray}
where $\omega(t_k)$ is a random variable. Both $\omega(t_k)$ and $z(t_k)$, especially the probability of survival of each agent at $t_k$,  define how the current index set $I(t_k)$ should be updated. 

As an example, let $\omega(t_k)$ be uniformly distributed on $[0,1]$ and  $J(t_k)\subseteq I(t_k)$ be the index set defined as
\begin{eqnarray}\label{J eqn}
  J(t_k) \triangleq &    \left\{ j \in I(t_k) \left| \  \omega(t_k) \geq \mbox{the probability of } \right. \right. \\
                           &     \left. \mbox{survival of agent} \ j \right\} \nonumber
\end{eqnarray}
Then the update $\psi$, can be defined as
\begin{equation}\label{update eqn}
 I(t_{k+1}) \ = \ \psi (I(t_k),z(t_k),\omega(t_k))\  \triangleq \ I(t_k) \setminus J(t_k) .
 \end{equation} 
Any changes on the index set $I(t_{k})$ will affect the entire swarm dynamics at the next time instance. 
Such coupling between locally deterministic dynamics (for each agent ${i}$) and globally stochastic dynamics (for the index set $I$ and time span $[t_0,t_f]$) is not well captured in standard deterministic or stochastic control. Furthermore, the performance metrics in adversarial swarms can  typically be expressed as a function of all agents at the final time $t_{f}$, i.e.,
\begin{eqnarray}
F(z (t_{f}))\label{eq:cost}
\end{eqnarray}
which is also stochastic as $z(t_{f})$ depends on the entire sequence of random variables $$\left\{ \omega(t_0),\omega(t_1),\cdots \omega(t_{f})\right\}.$$ This cost must thus be transformed to an expectation additionally dependent on the indicator set, $E\left[F(z (t_{f}), I(t_f))\right]$. We arrive at the following discrete stochastic optimal control problem 
\begin{eqnarray}
P0\triangleq\left\{
\begin{array}{l}
\min  \ \ J := E\left[F(z (t_{f}), I(t_f))\right] \\
 \textbf{subject to}  \nonumber \\
 \quad z(t_{k+1})  = \phi^kz(t_k), (t_k), I(t_k)), 
  \nonumber \\
 \quad I(t_{k+1}) = \psi (I(t_k),z(t_k),\omega(t_k)) \nonumber \\
 \quad  \mathcal{H}(z(t_k),u(t_k)) \le 0, 
\end{array} \right.
\label{OCP_A}
\end{eqnarray}
 that has not been well-addressed by the existing deterministic or stochastic optimal control frameworks. This optimal control problem has several distinctive features that make it challenging to solve:
\begin{itemize}
\item the random time-varying dimension of the dynamics and the state trajectory due to loss of agents in an adversarial environment;
\item probability dependent performance metric intertwining with locally deterministic agent dynamics (for each agent and at each local time instance), yet globally stochastic dynamical behavior (with respect to a time horizon and entire swarm);
\item high-dimension in both decision variables (e.g., trajectories of defending agents) and dynamical constraints (overall swarm dynamics that may include thousands of agents).
\end{itemize}

As an attempt to address these challenges, in the next section we propose some simplified models to approximate such optimal control problems.

\section{Proposed Solution Methods} \label{sec:soving-zombie}


In this section, we consider three alternative numerical schemes for solving the Problem $P0$.

\subsection{P1: Deterministic and Decoupled Optimal Control Problem Formulation}

We first consider the option of not modifying the swarm dynamics by the indicator set. This is the approach originally taken in~\cite{walton2018optimal}. This results in continuous dynamics over time. Rather than an indicator set, agent survival can also be modeled as continuous probabilities over time. Let $Q_j^I(t)$ by the probability of the $j$-th attacker surviving at time $t$, $Q_k^D(t)$  be the probability of the $k$-th defender surviving at time $t$, and $Q_0(t)$  be the probability of the HVU surviving. Survival probablities
\[
Q = [Q_0, Q_1^I, \dots Q_N^I, Q_1^D, \dots Q_M^D]
\]
can be modeled with the dynamics
\[
Q(t_{k+1}) = \Psi(Q(t_k), z(t_k))
\]
for example as the nonhomogeneous Poisson process of mutual attrition.

The continuous-time standard optimal control problem $P$ can be expressed in the following form
\begin{eqnarray} \label{OCP}
P1\triangleq\left\{
\begin{array}{l}
\min \ \ J := F(z(t_f), Q(t_f))    \\
 \textbf{subject to}   \\
 \quad z(t_{k+1}) = \phi (z(t_k),u(t_k))  \\
\quad Q(t_{k+1}) = \Psi(Q(t_k), z(t_k)) \\
 \quad  \mathcal{H}(z(t_k),u(t_k)) \le 0, 
\end{array} \right.
\end{eqnarray}
where $ F(z(t_f), Q(t_f))  $ represents the terminal cost, $\phi (z(t_k),u(t_k)$ represents the deterministic swarm dynamics, and $\mathcal{H}(z(t_k),u(t_k))$ is the constraints on the states and control inputs. 

This formulation, while diverging from the reality of the swarm scenario by ignoring attrition when considering dynamics, has the advantage of providing a smooth problem for optimization.

\subsection{P2: Weighted Forces Model}

Another approach, which we denote ``weighted force model,'' has the swarm dynamics depend on agent survival through the continuous time survival probabilities $Q(t)$ instead of the indicator set. This is done by weighting the contribution of each agent to the collective dynamics by its probability of survival. An example is provided in the simulation model description in equation \ref{eqn:weightedattackers}. Note that if survival probablities are all equal to $1$, this would return the dynamics of problem $P1$. If they were binary indicators of survival, this would return the dynamics of the indicator set coupled dynamics of $P0$. This approach maintains the smoothness properties of $P0$ while potentially better approximating $P0$. This model includes some unphysical characteristics (e.g., intra-swarm cohesion decreases as  agent survival probabilities decrease), but it has the essential characteristic that dead agents are no longer able to affect the dynamics of others. We summarize this as
\begin{eqnarray} \label{P2}
P2\triangleq\left\{
\begin{array}{l}
\min \ \ J := F(z(t_f), Q(t_f))    \\
 \textbf{subject to}   \\
 \quad z(t_{k+1}) = \phi (z(t_k),u(t_k), Q(t_k))  \\
\quad Q(t_{k+1}) = \Psi(Q(t_k), z(t_k)) \\
 \quad  \mathcal{H}(z(t_k),u(t_k)) \le 0, 
\end{array} \right.
\end{eqnarray}

\subsection{P3: Threshold model}

The threshold model treats all agents as fully alive until their survival probability drops below some threshold, which we choose as 50\%, after which they do not interact with other agents (in dynamics or attrition). This provides an update rule for equation (\ref{eq:I_dyn}), updating the indicator set $I(t_k)$. This update rule, however, uses the continuous probability dynamics. This provides a smooth problem in between changes to $I(t_{k+1})$. Furthermore, the index set is dependent on $Q$ rather than the random variable $\omega$, keeping the expectation out of the cost function $J$. Therefore, this problem also can be studied with gradient-based optimization schemes.
\begin{eqnarray} \label{P2}
P3\triangleq\left\{
\begin{array}{l}
\min \ \ J := F(z(t_f), I(t_f))    \\
 \textbf{subject to}   \\
 \quad z(t_{k+1}) = \phi^k (z(t_k),u(t_k), I(t_k))  \\
\quad Q(t_{k+1}) = \Psi(Q(t_k), z(t_k)) \\
 \quad I(t_{k+1}) = \psi (I(t_k), z(t_k), Q(t_k)) \nonumber \\
 \quad  \mathcal{H}(z(t_k),u(t_k)) \le 0, 
\end{array} \right.
\end{eqnarray}

\subsection{Validation: Monte Carlo}

To examine the efficacy of the three models, we test against Monte  Carlo simulations based on the defender trajectories generated by each optimization. These Monte Carlo simulations enact the agent and index set dynamics of $P0$ for the calculated defender controls. For the index set update of equation (\ref{eq:I_dyn}), we model the random variable $\omega(t_k)$ as a vector with an independent component for each agent including the HVU
\[
\omega = [\omega_0, \omega_1^I, \dots \omega_N^I, \omega_1^D, \dots \omega_M^D]
\]. 
This random variable is sampled at each $t_k$. Each agent component of the random variable is sampled uniformly.  The index set is updated by removing a set in the form of equation (\ref{J eqn}) with the definition
\begin{eqnarray}\label{J eqn}
  J(t_k) \triangleq &    \left\{ j \in I(t_k) \left| \  \omega_j(t_k) <\frac{Q_j(t_{k+1})}{Q_j(t_k)} \right. \right\} 
\end{eqnarray}
This ratio $\frac{Q_j(t_{k+1})}{Q_j(t_k)}$ models the probability in time interval $[t_k, t_{k+1}]$ that the $j$-th agent is destroyed.

\section{Case study: Defense Against a Swarm Attack}\label{sec:model}

To test the numerical schemes presented in the previous section, we consider a scenario where a swarm is attacking a high-value unit (HVU). The HVU is defended by a number of defending agents whose trajectories are described by a finite sum of Bernstein polynomials and are optimized to maximize the probability of the HVU survival. All agents are equipped with model weapons, such that all agents and the HVU have an attrition rate at each time step in the simulation that is determined by its relative position to all enemy agents. 


\subsection{Attacker Equations of Motion}
We consider $N$ attacking agents and $M$ defending agents, where attacking agent $i$ has position $x_i(t) \in R^3$ and defending agent $k$ has position $s_k(t) \in R^3$. The equation of motion for attacker $i$, where the acceleration $\ddot{x}_i$ at each time step is the sum of four forces, is
\begin{align} \label{eqn:attackers}
\ddot{x}_i =& \sum\limits_{j \ne i}^{{N}} {\frac{{{f_I}({x_{ij}})}}{{\left\| {{x_{ij}}} \right\|}}{x_{ij}}}  + \sum\limits_{k=1}^M {\frac{{{f_d}({s_{ik}})}}{{\left\| {{s_{ik}}} \right\|}}{s_{ik}}} \nonumber \\ 
& +  K\frac{h_i}{\left\|h_i\right\|} - {b}{\dot x_i},  \hspace{0.5in} i = 1\dots N
\end{align}
There are four terms in this equation, representing: (1) attractive and repulsive forces $f_I(x_{ij})$ from other attacking agents $j$, where $x_{ij} = x_i - x_j$ is the distance between attackers $i$ and $j$; (2) a constant ``virtual leader'' force with magnitude $K$ pulling them toward the HVU's position, where $h_i = h-x_i$ and $h$ is the position of the HVU; (3) purely repulsive forces $f_d(s_{ik})$ due to defending agents, where $s_{ik} = x_i - s_k$ is the distance between attacker $i$ and defender $k$; and (4) a damping force proportional to the $\dot{x}_i$.

For the mathematical forms of $f_I$ and $f_d$, we use a common model proposed by Leonard and Fiorelli~\cite{leonard2001virtual}, where $f_I$ and $f_d$ can be written as gradients of scalar potential functions that depend only on $x_{ij}$ and $s_{ik}$, respectively. Both terms include repulsive collision avoidance at short ranges, and $f_I$ includes attractive forces at intermediate ranges for swarm cohesion. Specifically, $f_I$ is repulsive when $\left\| {{x_{ij}}} \right\| \le {d_0}$, attractive when ${d_0} < \left\| {{x_{ij}}} \right\| \le {d_1}$, and zero when $\left\| {{x_{ij}}} \right\| > {d_1}$. Similarly, $f_d$ is repulsive when $\left\| {{s_{ik}}} \right\| \le {s_0}$ and zero when $\left\| {{s_{ik}}} \right\| > {s_0}$. To test robustness, we also performed simulations using the Reynolds dynamics model~\cite{reynolds1987flocks} instead of the first term in Eq.~\eqref{eqn:attackers} with qualitatively similar results to those we show here.

For the Threshhold model and the Monte Carlo simulations, the swarm dynamics of equation \ref{eqn:attackers} are modified by performing the summations only over the intersection of the original indices and the index set at time $t$. The weighted forces model weights the summation terms by the respective survival probabilities, as:

\begin{align} \label{eqn:weightedattackers}
\ddot{x}_i =& \sum\limits_{j \ne i}^{{N}} Q_j^I {\frac{{{f_I}({x_{ij}})}}{{\left\| {{x_{ij}}} \right\|}}{x_{ij}}}  + \sum\limits_{k=1}^M Q_k^D {\frac{{{f_d}({s_{ik}})}}{{\left\| {{s_{ik}}} \right\|}}{s_{ik}}} \nonumber \\ 
& +  K\frac{h_i}{\left\|h_i\right\|} - {b}{\dot x_i},  \hspace{0.5in} i = 1\dots N
\end{align}

\subsection{Defender Equations of Motion}

In this paper we have used a 3D double integrator model to represent defender dynamics.
\begin{align} \label{defenders}
\ddot{s}_k = u_k, \hspace{0.5in} k = 1\dots M
\end{align}
where $u_k (t) \in R^3$ and absolute value of each element of $u_k$, $(|u_{kj}|, j =1,2,3)$  is bounded by $u_{\rm max}$. 

The discrete dynamics of defenders and attackers results from explicit discretization of these continuous dynamics.

\subsection{Mutual Attrition model}

To model mutual attrition between enemy agents, we choose a pairwise damage function that takes as an argument the relative position between the two agents. This function has a value of 1 when its argument is 0 (i.e., when the agents are at the same position), and the function smoothly and continuously decreases as the argument increases. We use an inverted cumulative normal distribution, which we denote $\Phi$, to accomplish this, but our results are highly insensitive to this choice. Thus, the rate at which attacker $i$ is destroyed due to defender $k$ is $d_{ik}^{\rm att} = {\lambda _d}{\Phi}({\|s_{ik}\|^2}/{\sigma _d})$, where $\sigma _d$ is a range parameter and $\lambda _d$ is a rate-of-fire parameter. Similarly, the attrition rate of defender $k$ due to attacker $i$ is $d_{ki}^{\rm def} = {\lambda _a}{\Phi}\left({{{{\|s_{ik}\|^2}}}}/{{{\sigma _a}}}\right)$, and the attrition rate of the HVU is $d_i^{\rm hvu} = \lambda_a \Phi\left({\|h_i\|^2}/{\sigma _a}\right)$, where $\sigma _a$ and $\lambda_a$ correspond to the range and rate of fire of the attackers' weapons.

At each time step, all attackers are firing at all defenders, and vice versa. Thus, for example, the probability that attacker $j$ would survive during time interval $\Delta t$ can be written as $\prod\limits_i^{{M_{}}} {(1 - \left[ {d_{ji}^{\rm att}Q_i^D (t)} \right]} \Delta t)$. The survival probabilities ${Q_j^I}(t+\Delta t)$ for attacker $j$, $Q_i^D(t)$ for defender $i$, and $Q_0(t)$ for the HVU are governed by 
\begin{align}\label{eq:Q_i(t)} \nonumber
{Q_j^I}(t_{k+1}) &= {Q_j^I}(t_k) {\prod\limits_i^{{M_{}}} {(1 -  {d_{ji}^{\rm att}Q_i^D (t_k)}} \left[ t_{k+1}-t_k  \right])}, \\
Q_i^D(t_{k+1}) &=Q_i^D(t_k){\prod\limits_j^{{N_{}}} {(1 - \left[ {d_{ij}^{\rm def}{Q_j^I}(t_k)} \right]} \left[ t_{k+1}-t_k \right])}, \nonumber \\
Q_0(t_{k+1}) &=Q_0(t_k){\prod\limits_j^{{N_{}}} {(1 - \left[ {d_{ij}^{\rm def}{Q_j^I}(t_k)} \right]} \left[ t_{k+1}-t_k  \right])},
\end{align}
Initial conditions are set to ${Q_i}(0) = 1$ for all agents and the HVU.

\subsection{Numerical Methods}

For a given set of defender trajectories, we numerically integrate Eqs.~\eqref{eqn:attackers} and \eqref{eq:Q_i(t)} using standard numerical methods of molecular dynamics (MD)~\cite{allen2017computer}. In particular, we integrate Eq.~\eqref{eqn:attackers} using a velocity-Verlet integration scheme~\cite{Verlet1967}. To solve problems $P1$, $P2$, and $P3$, we define the terminal cost as the probability that the HVU was destroyed at the end of the simulation. The direct methods of optimal control can then be applied to problems $P1$, $P2$, and $P3$ to find defender trajectories that minimize the probability that the HVU is destroyed. 

\begin{figure}
\raggedright 
(a) \\ \centering \includegraphics[trim=5mm 2mm 5mm 5mm, clip, width=0.8\columnwidth]{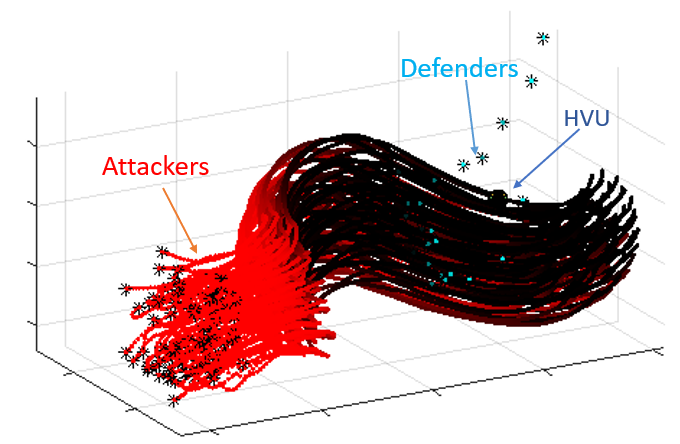} \\
\raggedright (b) \\ \centering \includegraphics[trim=2mm 2mm 5mm 5mm, clip, width=0.8\columnwidth]{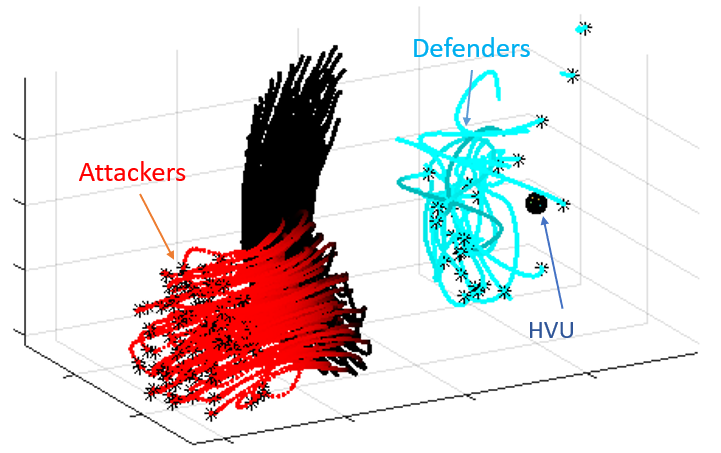} \\
\raggedright(c) \\ \centering \includegraphics[width=\columnwidth]{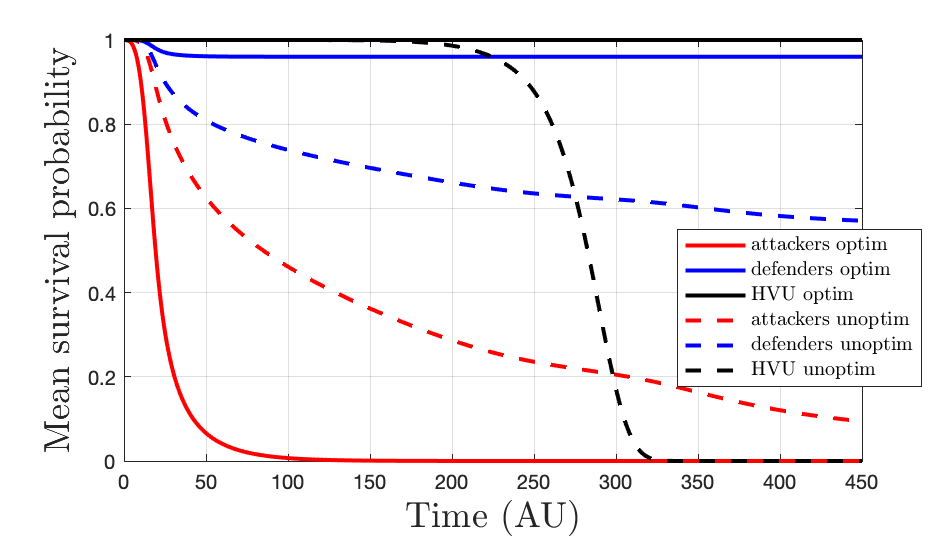}
\caption{Unoptimized (a) and optimized (b) defender trajectories are shown for a confrontation of 100 attacking swarm agents with a HVU protection force of 25 defenders with superior weapons. The optimized trajectories defend the HVU more effectively, as shown in (c).}
\label{fig:opt-unopt}
\end{figure}

Figure~\ref{fig:opt-unopt} shows an optimization of$P1$, where survival probability is not coupled to the equations of motion of the attackers, for 25 defenders protecting an HVU against a swarm of 100 attackers. The defenders have a 50\% larger weapons range ($\sigma_d/\sigma_a = 1.5$) as well as double the fire rate with respect to the attackers ($\lambda_d/\lambda_a = 2$). Figure~\ref{fig:opt-unopt}(a) shows the results of a simulation where the defenders remain in place, and Figure~\ref{fig:opt-unopt}(b) shows results of a simulation after optimizing the trajectories of the defenders. The attackers are red and the defenders are cyan, but turn to black as their survival probability decreases to zero. In both scenarios the defenders are stronger and suffer fewer losses; however, in the unoptimized scenario, some of the attacking agents manage to penetrate the defenders’ zone and destroy the HVU. Figure \ref{fig:opt-unopt}(c) shows that the optimized trajectories lead to better results for the survival of the HVU, the survival of defenders, and the defeat of the attackers (although the survival of the HVU is the only metric used in the optimization). We also note that the defending agents utilize both herding and weapons, as seen in Fig.~\ref{fig:opt-unopt}.


\subsection{Comparing  Performance of the Proposed Models}\label{sec:optimization comparizon}


Similar results to those shown in Fig.~\ref{fig:opt-unopt} are found when considering problems $P2$ (where spatial interactions are multiplied by the survival probabilities) and $P3$ (where spatial interactions and damage functions are turned off below a 50\% survival probability threshold). These results are qualitatively indistinguishable from the results we show for $P1$, so we do not show them here. However, a key question remains: how well do the modeling frameworks used in problems $P1$ (uncoupled dynamics), $P2$, and $P3$ compare with the stochastic problem $P0$? Problems $P1$, $P2$, and $P3$ are approximations to $P0$. 

To this end, we compare results of a single simulation with fixed defender trajectories using the modeling frameworks corresponding to each problem. We consider a case with $M=50$ defenders and $N=50$ attackers with identical rates of fire ($\lambda_a = \lambda_d$). Attackers, however, have a 10\% larger range, $\sigma_a/\sigma_d=1.1$. The defender trajectories are obtained by optimization of problem $P2$, and these trajectories are then fed into each of the four types of modeling framework. Results for $P0$ (Monte Carlo) are averaged over 200 simulations; all other results use a single simulation, since there is no randomness.

\begin{figure*}
\centering
\includegraphics[trim=40mm 0mm 40mm 0mm, clip, width=\textwidth]{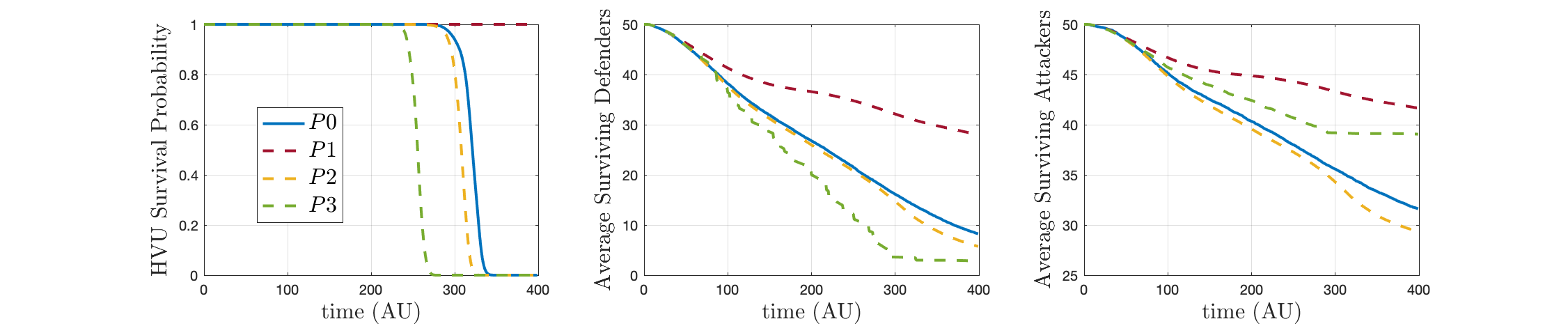}
\caption{Comparison of the performance of each simulation method using fixed defender trajectories. The simulation methods for $P2$ and $P3$ that couple attrition and spatial dynamics agree much better with the random (Monte Carlo) simulations for $P0$.}
\label{fig:compare-Ps}
\end{figure*}

Figure~\ref{fig:compare-Ps} shows results from each type of simulation, specifically the mean survival probabilities of attackers, defenders, and the HVU. This figure shows that $P2$ and $P3$ have similar results to $P0$, which is typical of all simulations. However, $P1$ does not agree with $P0$, $P2$, or $P3$. Instead, $P1$ greatly overestimates the probability of HVU survival, which is also typical of all simulations. Physically, this has an obvious explanation: $P1$ had no coupling between attrition and spatial dynamics, so attackers would still try to avoid defenders who were in their path, even if these defenders had a very low probability of survival. In contrast, the modeling frameworks corresponding to $P2$ and $P3$ reduced the spatial interactions as survival probability decreased. So, even with the simple coupling used in $P2$ and $P3$ and essentially no calibration, these approximations agreed relatively well with the stochastic results in problem $P0$. This result highlights our main point in this paper, which is that modeling and control frameworks for adversarial autonomy must include attrition, and the attrition modeling should be coupled to the spatial dynamics of the agents.

\section{Conclusions}\label{sec:conclusions}

In this paper we have addressed the question of optimal motion planning for large-scale autonomous systems that explicitly include attrition. We proposed a novel optimal control problem that explicitly includes random reduction of an index set of surviving agents in time. Since this problem is inherently stochastic, we proposed three ``smooth'' approximations that can be solved using direct methods of optimal control. By considering a case study of defending agents protecting an HVU from an attacking swarm, we showed that these approximations can be solved and that they give results that are consistent with the stochastic problem, especially if the attrition and spatial dynamics are coupled. We note that our results assume that the cooperating strategies and weapons capabilities of the attacking swarm are known or can be estimated. Estimation can be considered seperately, using an approach such as in \cite{JGCD2019}. Parameter uncertainty can also be added into this framework, using methods such as in \cite{walton_IJC}. The framework we describe here can be applied to an entire class of adversarial autonomy problems. For example, attrition could result from many factors, including environmental or terrain-related causes. Future work might focus on improving the approximation methods ($P1$, $P2$, and $P3$) as well as improving stochastic optimization techniques such that $P0$ could be treated directly.
\section{ACKNOWLEDGMENTS}
This work was supported in part by ONR SoA program and by NPS CRUSER program. 
\bibliographystyle{IEEEtran}
\bibliography{references}  

\begin{thebibliography}{10}
\providecommand{\url}[1]{#1}
\csname url@samestyle\endcsname
\providecommand{\newblock}{\relax}
\providecommand{\bibinfo}[2]{#2}
\providecommand{\BIBentrySTDinterwordspacing}{\spaceskip=0pt\relax}
\providecommand{\BIBentryALTinterwordstretchfactor}{4}
\providecommand{\BIBentryALTinterwordspacing}{\spaceskip=\fontdimen2\font plus
\BIBentryALTinterwordstretchfactor\fontdimen3\font minus
  \fontdimen4\font\relax}
\providecommand{\BIBforeignlanguage}[2]{{%
\expandafter\ifx\csname l@#1\endcsname\relax
\typeout{** WARNING: IEEEtran.bst: No hyphenation pattern has been}%
\typeout{** loaded for the language `#1'. Using the pattern for}%
\typeout{** the default language instead.}%
\else
\language=\csname l@#1\endcsname
\fi
#2}}
\providecommand{\BIBdecl}{\relax}
\BIBdecl

\bibitem{evers2014robust}
L.~Evers, T.~Dollevoet, A.~I. Barros, and H.~Monsuur, ``Robust uav mission
  planning,'' \emph{Annals of Operations Research}, vol. 222, no.~1, pp.
  293--315, 2014.

\bibitem{reynolds1987flocks}
C.~W. Reynolds, \emph{Flocks, herds and schools: A distributed behavioral
  model}.\hskip 1em plus 0.5em minus 0.4em\relax ACM, 1987, vol.~21.

\bibitem{vicsek}
\BIBentryALTinterwordspacing
T.~Vicsek, A.~Czir\'ok, E.~Ben-Jacob, I.~Cohen, and O.~Shochet, ``Novel type of
  phase transition in a system of self-driven particles,'' \emph{Phys. Rev.
  Lett.}, vol.~75, pp. 1226--1229, Aug 1995. [Online]. Available:
  \url{https://link.aps.org/doi/10.1103/PhysRevLett.75.1226}
\BIBentrySTDinterwordspacing

\bibitem{leonard2001virtual}
N.~E. Leonard and E.~Fiorelli, ``Virtual leaders, artificial potentials and
  coordinated control of groups,'' in \emph{Proceedings of the 40th IEEE
  Conference on Decision and Control (Cat. No. 01CH37228)}, vol.~3.\hskip 1em
  plus 0.5em minus 0.4em\relax IEEE, 2001, pp. 2968--2973.

\bibitem{brambilla2013swarm}
M.~Brambilla, E.~Ferrante, M.~Birattari, and M.~Dorigo, ``Swarm robotics: a
  review from the swarm engineering perspective,'' \emph{Swarm Intelligence},
  vol.~7, no.~1, pp. 1--41, 2013.

\bibitem{TR2018}
S.-J. Chung, A.~A. Paranjape, P.~Dames, S.~Shen, and V.~Kumar, ``A survey on
  aerial swarm robotics,'' \emph{IEEE Transactions on Robotics}, vol.~34,
  no.~4, pp. 837--855, August 2018.

\bibitem{herdingbirds}
A.~A. Paranjape, S.-J. Chung, K.~Kim, and D.~H. Shim, ``Robotic herding of a
  flock of birds using an unmanned aerial vehicle,'' \emph{IEEE Transactions on
  Robotics}, vol.~34, no.~4, pp. 901--915, 2018.

\bibitem{Egerstedt}
M.~A. Haque, A.~R. Rahmani, and M.~B. Egerstedt, ``A hybrid, multi-agent model
  of foraging bottlenose dolphins,'' in \emph{IFAC Proceedings Volumes},
  vol.~42, 2009, pp. 262--267.

\bibitem{sheepdogs1}
D.~Str{\"o}mbom, R.~P. Mann, A.~M. Wilson, S.~Hailes, A.~J. Morton, D.~J.
  Sumpter, and A.~J. King, ``Solving the shepherding problem: heuristics for
  herding autonomous, interacting agents,'' \emph{Journal of the royal society
  interface}, vol.~11, no. 100, p. 20140719, 2014.

\bibitem{sheepdogs2}
A.~Pierson and M.~Schwager, ``Bio-inspired non-cooperative multi-robot
  herding,'' in \emph{2015 IEEE International Conference on Robotics and
  Automation (ICRA)}.\hskip 1em plus 0.5em minus 0.4em\relax IEEE, 2015, pp.
  1843--1849.

\bibitem{walton2018optimal}
C.~Walton, P.~Lambrianides, I.~Kaminer, J.~Royset, and Q.~Gong, ``Optimal
  motion planning in rapid-fire combat situations with attacker uncertainty,''
  \emph{Naval Research Logistics (NRL)}, vol.~65, no.~2, pp. 101--119, 2018.

\bibitem{allen2017computer}
M.~P. Allen and D.~J. Tildesley, \emph{Computer simulation of liquids (2nd
  Ed.)}.\hskip 1em plus 0.5em minus 0.4em\relax Oxford university press, 2017.

\bibitem{Verlet1967}
\BIBentryALTinterwordspacing
L.~Verlet, ``Computer "experiments" on classical fluids. i. thermodynamical
  properties of lennard-jones molecules,'' \emph{Phys. Rev.}, vol. 159, pp.
  98--103, Jul 1967. [Online]. Available:
  \url{https://link.aps.org/doi/10.1103/PhysRev.159.98}
\BIBentrySTDinterwordspacing

\bibitem{JGCD2019}
Q.~Gong, W.~Kang, C.~Walton, I.~Kaminer, and H.~Park, ``Partial observability
  analysis of an adversarial swarm model,'' \emph{Journal of Guidance, Control,
  and Dynamics}, pp. 1--12, 2019.

\bibitem{walton_IJC}
C.~Walton, I.~Kaminer, and Q.~Gong, ``Consistent numerical methods for state
  and control constrained trajectory optimisation with parameter dependency,''
  \emph{International Journal of Control}, 2020.

\end{thebibliography}
\end{document}